\nonstopmode
\documentclass[10pt,reqno]{amsart}
\usepackage{latexsym}
\usepackage{fancyhdr}
\usepackage{amsmath, amssymb}
\usepackage[ansinew]{inputenc}
\usepackage[all]{xy}
\usepackage{pdflscape}
\usepackage{longtable}
\usepackage{rotating}
\usepackage{verbatim}
\usepackage{hyperref}
\usepackage{subfigure}
\usepackage{mathrsfs}
\usepackage{tensor}
\usepackage{enumerate}
\usepackage{mdwlist}

\swapnumbers
\theoremstyle{plain}
\newtheorem*{lemma*}{Lemma}
\newtheorem{lemma}[subsection]{Lemma}
\newtheorem*{theorem*}{Theorem}
\newtheorem{theorem}[subsection]{Theorem}
\newtheorem*{proposition*}{Proposition}
\newtheorem{proposition}[subsection]{Proposition}
\newtheorem*{corollary*}{Corollary}
\newtheorem{corollary}[subsection]{Corollary}
\newtheorem*{claim*}{Claim}
\theoremstyle{definition}
\newtheorem*{definition*}{Definition}

\newtheorem*{example*}{Example}
\newtheorem{example}[subsection]{Example}

\newtheorem{algorithm}[subsection]{Algorithm}
\newtheorem*{algorithm*}{Algorithm}
\newtheorem*{remark*}{Remark}
\newtheorem{remark}[subsection]{Remark}
\newtheorem{remarks}[subsection]{Remarks}

\makeatletter
  \let\c@equation\c@subsection
  
\makeatother

\newenvironment{demo}[1]{\par\smallskip\noindent{\bf #1.}}{\par\smallskip}

\numberwithin{equation}{section}

\def\al{\alpha}
\def\be{\beta}
\def\ga{\gamma}

\def\ep{\epsilon}

\def\la{\lambda}

\def\si{\sigma}

\def\ta{\tau}

\def\om{\omega}

\def\Ga{\Gamma}
\def\De{\Delta}
\def\Th{\Theta}
\def\La{\Lambda}

\def\Ph{\Phi}

\def\C{\mathbb{C}}

\def\I{\mathbb{I}}
\def\N{\mathbb{N}}
\def\Q{\mathbb{Q}}
\def\R{\mathbb{R}}
\def\Z{\mathbb{Z}}

\def\p{\partial}

\def\<{\langle}
\def\>{\rangle}
\renewcommand{\o}{\circ}

\let\on=\operatorname
\newcommand{\sr}[1]%
{\ifmmode{}^\dagger\else${}^\dagger$\fi\ifvmode
\vbox to 0pt{\vss
 \hbox to 0pt{\hskip\hsize\hskip1em
 \vbox{\hsize3cm\raggedright\pretolerance10000
 \noindent #1\hfill}\hss}\vss}\else
 \vadjust{\vbox to0pt{\vss%
 \hbox to 0pt{\hskip\hsize\hskip1em%
 \vbox{\hsize3cm\raggedright\pretolerance10000%
 \noindent #1\hfill}\hss}\vss}}\fi%
}

\newcommand{\tC}[2]{\tensor[^{#1}]{C}{^{#2}}}

\title[Differentiable roots, eigenvalues, and eigenvectors]
{Differentiable roots, eigenvalues, and eigenvectors}

\author[A. Rainer]
{Armin Rainer}

\address{Armin Rainer: Fakult\"at f\"ur Mathematik, Universit\"at Wien, 
Nordbergstrasse~15, A-1090 Wien, Austria}

\email{armin.rainer@univie.ac.at}

\begin{document}

\begin{abstract} 
  We determine the conditions for the existence of $C^p$-roots of curves of monic complex polynomials
  as well as for the existence of $C^p$-eigenvalues and $C^p$-eigenvectors of curves of normal complex matrices.  
\end{abstract}


\thanks{The author was supported by the Austrian Science Fund (FWF), Grant P~22218-N13}
\date{January 22, 2013}

\maketitle

\section{Introduction}

Consider a monic polynomial $P$ whose coefficients are complex valued $C^p$-functions, 
where $p \in \N \cup \{\infty\}$, 
defined in an 
open interval $I \subseteq \R$:
\begin{equation} \label{eq:P}
	P(t)(z) = z^n+\sum_{j=1}^n (-1)^j a_j(t) z^{n-j}, \quad a_j \in C^p(I,\C) \text{ for all } j.
\end{equation}
We shall say that $P$ is a \emph{$C^p$-curve of polynomials}. Furthermore, we shall say that $P$ is 
\emph{$C^q$-solvable}, 
if there exists a $C^q$-parameterization of its roots, i.e., there are functions $\la_j \in C^q(I,\C)$, $1 \le j \le n$, 
so that 
\begin{equation*} 
	P(t)(z) = \prod_{j=1}^n (z-\la_j(t)), \quad t \in I.
\end{equation*}
In general a $C^\infty$ or even a real analytic curve of polynomials need not have differentiable roots, for instance
$P(t)(z) = z^2-t$. In this example the order of contact of the roots is too low. 
On the other hand ``low'' order of contact of the roots is excluded
if $P$ is \emph{hyperbolic} (that is all roots are real). However,  
it is well-known that a general $C^\infty$-curve of hyperbolic polynomials is not $C^{1,\al}$-solvable for any $\al>0$, 
see \cite{BBCP06};
the loss of smoothness is caused by oscillation and infinite (i.e.\ ``high'') order of contact of the roots. 
On the other hand $C^n$ (resp.\ $C^{2n}$)-coefficients guarantee $C^1$ (resp.\ twice differentiable) solvability,
see \cite{ColombiniOrruPernazza12}. 
A sufficient condition for $C^\infty$-solvability of hyperbolic polynomials is that no two distinct roots have infinite 
order of contact, see \cite{AKLM98},
another sufficient condition is definability of the coefficients which allows for infinite contact but excludes 
oscillation, see \cite{RainerOmin}.  
We use the notion \emph{definability} always with respect to a fixed arbitrary o-minimal expansion of the real field;
cf.\ the book \cite{vandenDries98} or the expository paper \cite{vandenDriesMiller96}. 
See \cite[Table 1]{RainerQA} for an overview of the state of knowledge on the problem of choosing the roots of polynomials in a smooth way. 

E.\ Bierstone suggested to find
the conditions for $C^q$-solvability, for arbitrary $q$, 
in terms of the differentiability of the coefficients and the (finite) order of contact of the roots.
We solved this problem in \cite{RainerOmin} for hyperbolic polynomials $P$ whose coefficients 
are definable and for very special ``non-definable'' hyperbolic polynomials.
In this paper we give a complete solution without assuming hyperbolicity or definability, see Theorem~\ref{main}.

To this end we associate with any germ at $t \in I$ of a $C^p$-curve of polynomials $P$ as in \eqref{eq:P} 
a labeled rooted tree $T(P,t)$ which reflects the iterated factorization 
of $P$ near $t$, see Algorithm~\ref{alg}. 
In order to avoid pathological situations we restrict our attention to so-called \emph{admissible} trees $T(P,t)$, see \ref{gbf}. 
We say that $T(P,t)$ is \emph{good} if this algorithm produces a 
a splitting into linear factors; this is the case if the order of contact of the roots is sufficiently high at each 
step, e.g., if $P$ is hyperbolic. Otherwise $T(P,t)$ is called \emph{bad}. Some bad rooted trees we call \emph{fair} 
if they have certain properties that are convenient for our goal, see \ref{gbf}. 
Via the trees $T(P,t)$ we associate with $P$ numbers 
$\Ga(P), \hat \Ga(P), \ga(P) \in \N\cup \{\infty\}$ which encode the conditions for solvability.
More precisely, our main Theorem~\ref{main} states that under the assumption that roots do not have infinite order 
of contact:
\begin{itemize}
   \item For all $p \in \N \cup \{\infty\}$, 
   $P$ is $C^{p+\ga(P)}$-solvable, if $P$ is $C^{p+\Ga(P)}$ and $T(P,t)$ is good for all $t$,
   \item $P$ is $C^{\ga(P)}$-solvable, if $P$ is $C^{\Ga(P)}$ and $T(P,t)$ is fair for all $t$, 
   \item $P$ is $C^{\ga(P)}$-solvable, if $P$ is $C^{1+\hat \Ga(P)}$ and $T(P,t)$ is bad for some $t$.
\end{itemize}
In Section~\ref{sec:ex} we give examples which show that the conditions in this statement are optimal, except perhaps in the 
third item. 
As a corollary we recover the results of \cite{Reichard80} on radicals of functions, 
see Corollary~\ref{cor1}. 

In Section~\ref{sec:normal} we solve the corresponding problem for curves of normal complex matrices 
$A(t)=(A_{ij}(t))_{1\le i,j\le n}$, $t \in I$, whose eigenvalues have finite order of contact. 
Here we define a number $\Th(A) \in \N \cup \{\infty\}$ that satisfies $\Th(A) \le \ga(P_A) \le \Ga(P_A)$,
where $P_A$ is the characteristic polynomial of $A$, 
and we prove in Theorem~\ref{normal} (under a similar admissibility condition) that: 
\begin{itemize}
  \item For all $p \in \N \cup \{\infty\}$, 
  the eigenvalues of $A$ can be parameterized by $C^{p+\Th(P)}$-functions, the eigenvectors 
  by $C^{p}$-functions, if $A$ is $C^{p+\Th(P)}$.   
\end{itemize}
These conditions are optimal, see Example~\ref{ex:normal}.
We want to stress the fact that here no loss of smoothness occurs and, loosely speaking, all admissible trees $T(A,t)$
are good. Similar phenomena have been observed in \cite{RainerN}.

In Section~\ref{sec:def} and Section~\ref{sec:defnormal} we state versions for definable 
polynomials and definable normal matrices, 
where admissibility of the associated trees and the assumption on the finite order of contact of the roots or eigenvalues are not needed.
This works well for the eigenvalues but not for the eigenvectors, which generally do not admit continuous
parameterizations without these assumptions, see Example~\ref{ex:normal2}.  

Let us explain by means of a very simple example the principles behind these results and its proofs.
Consider the polynomial $P(t)(z) = z^2-f(t)$, where $f$ is a germ at $0 \in \R$ of a $C^{m+2p}$-function so that 
$f(t)=t^{2p} g(t)$ and $g(0) \ne 0$. The roots of $P$ are given by $\pm t^p \sqrt{g(t)}$.   
The tree $T(P,0)$, that consists of the root with label $(2,p)$ and two leaves with label $(1,0)$, is good. 
We have $\Ga_0(P)=2p$ and $\ga_0(P) = p$. 
Obviously, $g$ is $C^m$ and, as $g(0) \ne 0$, so is $\sqrt{g}$. But $g$ is also what we call a 
\emph{$\tC{2p}{m+2p}$-function}, i.e., 
a function defined near $0 \in \R$ that becomes $C^{m+2p}$ if multiplied by the monomial
$t^{2p}$. In particular, $g$ is a $\tC{p}{m+p}$-function and thus also $\sqrt{g}$, since  
for any $\tC{p}{m+p}$-function $g$ and any $C^{m+p}$-germ $F$ at $g(0)$ the composite $F \o g$ is a 
$\tC{p}{m+p}$-function.
So the roots $\pm t^p \sqrt{g(t)}$ are $C^{m+p}$.  
We introduce and discuss $\tC{p}{m}$-functions in Section~\ref{sec:pCm}.

\subsection*{Notation and conventions}
We use $\N = \N_{>0} \cup \{0\}$. By convention $\infty + \Z=\infty$.
For $r \in \R$, we denote by $\lfloor r \rfloor$ (resp.\ $\lceil r \rceil$) the largest (resp.\ smallest) $p \in \Z$ such that 
$p \le r$ (resp.\ $p \ge r$).
By $\{r\} := r-\lfloor r \rfloor$ we mean the fractional part of $r$. 
The identity matrix of size $n$ is denoted by $\I_n$.

For a continuous complex valued function $f$ defined near $t_0$ in $\R$, 
let the \emph{multiplicity} $m_{t_0}(f)$ at $t_0$ be the supremum 
of all integers $p$ such that $f(t)=(t-t_0)^p g(t)$ near $t_0$ for a continuous function $g$. 
Note that, if $f$ is of class $C^n$ and $m_{t_0}(f) < n$, then $f(t) = (t-t_0)^{m_{t_0}(f)} g(t)$ 
near $t_0$, where $g$ is $C^{n-m_{t_0}(f)}$ and $g(t_0) \ne 0$. 

For a tuple $a=(a_1,\ldots,a_n)$ of germs at $0 \in \R$ of complex valued $C^0$-functions satisfying 
$m_0(a_k) \ge kr$ for some $r \in \N$, we say that $a$ is \emph{$r$-divisible} and define $a_{(r)}$ by setting
\[
  a_{(r)}(t):= (t^{-r}a_1(t),t^{-2r}a_2(t),\ldots,t^{-2n}a_n(t)).
\]
If $a_{(r)}(0) \ne 0$ we say that $a$ is \emph{strictly $r$-divisible}.

A monic polynomial $P$ of degree $n$ will frequently be identified with the tuple 
$a(P) = (a_1(P),\ldots,a_n(P))$ of its coefficients 
so that $P$ takes the form \eqref{eq:P} with $a_j=a_j(P)$. 

A rooted tree is a tree with a fixed special vertex, the root.
Writing $V \le W$ if $V$ belongs to the path between $W$ and the root defines a natural partial order on 
the set of vertices. The successors of a vertex $V$ are all vertices $W \ge V$ connected to $V$ by an edge. 
The maximal elements are called leaves. 
The height of a vertex in a rooted tree is the 
number of edges in the path that connects the vertex to the root. 
A rooted tree is called trivial if it consists just of its root.

\section{Preliminaries on polynomials}

For a monic polynomial with complex coefficients $a_1,\ldots,a_n$ 
and roots $\la_1,\ldots,\la_n$, 
\begin{equation*} \label{P}
P(z) = z^n + \sum_{j=1}^n (-1)^j a_j z^{n-j} = \prod_{j=1}^n (z-\la_j),
\end{equation*}
we have 
\[
 a_i = \sigma_i(\la_1,\ldots,\la_n)	= \sum_{1 \le j_1 < \cdots < j_i \le n} \la_{j_1} \cdots \la_{j_i}
\]
by Vieta's formulas, where $\sigma_1,\ldots,\sigma_n$ 
are the elementary symmetric functions in $n$ variables.
The Newton polynomials 
$s_i = \sum_{j=1}^n \la_j^i$, $i \in \N$,  
form a different set of generators for the algebra of symmetric polynomials on $\C^n$ and satisfy
\begin{equation*} \label{rec}
s_k - s_{k-1} \sigma_1 + s_{k-2} \sigma_2 - \cdots 
+ (-1)^{k-1} s_1 \sigma_{k-1} + (-1)^k k \sigma_k = 0, \quad (k \ge 1).
\end{equation*}
Consider the so-called Bezoutiant
\[ 
B := 
\begin{pmatrix} 
s_0 & s_1 & \ldots & s_{n-1}\\ 
s_1 & s_2 & \ldots & s_n \\ 
\vdots & \vdots & \ddots & \vdots\\ 
s_{n-1} & s_n & \ldots &  s_{2n-2} 
\end{pmatrix} 
= \left(s_{i+j-2}\right)_{1 \le i,j \le n}.
\]
If $B_k$ denotes the minor formed by the first $k$ rows and columns of $B$, then we have 
\begin{equation} \label{eq:del}
\Delta_k(\la) := \on{det} B_k(\la) = \sum_{i_1 < \cdots < i_k} 
(\la_{i_1}-\la_{i_2})^2 \cdots (\la_{i_1}-\la_{i_k})^2 \cdots (\la_{i_{k-1}}-\la_{i_k})^2
\end{equation}
and there exist unique polynomials $\tilde{\Delta}_k$ satisfying
$\Delta_k = \tilde{\Delta}_k \circ (\si_1,\ldots,\si_n)$. 
The number of
distinct roots of $P$ equals the
maximal $k$ such that $\tilde \De_k(P) \ne 0$.

\begin{lemma}[Splitting lemma {\cite[3.4]{AKLM98}}]  \label{split}
Let $P_0$ be a monic polynomial satisfying 
$P_0 = P_1 \cdot P_2$, where $P_1$ and $P_2$ are polynomials without common root. 
Then for $P$ near $P_0$ we have $P = P_1(P) \cdot P_2(P)$ 
for analytic mappings 
of monic polynomials $P \mapsto P_1(P)$ and $P \mapsto P_2(P)$, 
defined for $P$ 
near $P_0$, with the given initial values.
\end{lemma}

\begin{lemma}[{\cite[II Thm.\ 5.2]{Kato76}}] \label{controots}
  Any $C^0$-curve of polynomials as in \eqref{eq:P} is $C^0$-solvable.
\end{lemma}

\begin{lemma} \label{glue}
Let $P(t)$, $t \in I$, be a $C^0$-curve of polynomials as in \eqref{eq:P} and	
let $p \in \N$. 
If $P$ is locally $C^p$-solvable, then $P$ is globally $C^p$-solvable. 
\end{lemma}

\begin{demo}{Proof}
Let $\la=(\la_1,\ldots,\la_n) : I \supsetneq (a,b) \to \C^n$ be a maximal $C^p$-parameterization of the roots so that $b \in I$.
By assumption there is a local $C^p$-parameterization $\mu=(\mu_1\ldots,\mu_n)$ 
of the roots near $b$.
Let $t_0$ be in the common domain of $\la$ and $\mu$ and consider a sequence $t_k \to t_0^-$.
For each $k$ the $n$-tuples $\la(t_k)$ and $\mu(t_k)$ differ just by a permutation, and by passing to a subsequence we can assume 
that this permutation $\ta$ is independent of $k$, i.e., $\la(t_k) = \ta.\mu(t_k)$ for all $k$. 
Rolle's theorem implies that $\la^{(q)}(t_0) = \ta.\mu^{(q)}(t_0)$ for all $0 \le q \le p$, and hence  
$\la$ has a $C^p$-extension beyond $b$, contradicting its maximality.
\qed\end{demo}

\section{\texorpdfstring{$\tC{p}{m}$}{pCm}-functions} \label{sec:pCm}

\subsection{\texorpdfstring{$\tC{p}{m}$}{pCm}-functions}
Let $p,m \in \N$ with $p \le m$.
A continuous complex valued function $f$ defined near $0 \in \R$ is called a \emph{$\tC{p}{m}$-function} 
if $t \mapsto t^p f(t)$ belongs to $C^m$.
A vector or matrix valued function $f=(f_i)_{1\le i \le n}$ is called $\tC{p}{m}$ 
if each component $f_i$ is $\tC{p}{m}$.
(We will not use the notation $f \in \tC{p}{m}$ if $p>m$ or $f \not\in C^0$.)

Let $I \subseteq \R$ be an open interval containing $0$. 
Then $f : I \to \C$ is $\tC{p}{m}$ if and only if 
it has the following properties, cf.\ \cite[4.1]{Spallek72}, \cite[Satz~3]{Reichard79}, or \cite[Thm.~4]{Reichard80}:
\begin{itemize}
\item $f \in C^{m-p}(I)$.
\item $f|_{I\setminus \{0\}} \in C^m(I\setminus \{0\})$.
\item $\lim_{t \to 0} t^k f^{(m-p+k)}(t)$ exists as a finite number for all $0 \le k \le p$. 
\end{itemize}

\begin{proposition} \label{prop}
Let $k,l,m,p,q \in \N$.
We have:
\begin{enumerate}
\item[\thetag{1}] If $f$ is $\tC{p}{m}$, then $f$ is also $\tC{p-k}{m-k-l}$ for all $0 \le k \le p$, $0 \le l \le m-p$.
\item[\thetag{2}] If $f$ is $\tC{p}{m+p+q}$, then $f$ is also $\tC{p+q}{m+p+q}$.
\item[\thetag{3}] If $g=(g_1,\ldots,g_n)$ is $\tC{p}{m}$ and $F$ is $C^m$ near $g(0)\in \C^n$, then $F \o g$ is $\tC{p}{m}$.
\item[\thetag{4}] If $g_1,\ldots,g_n$ are $\tC{p}{m}$, then $\sum g_i$ and $\prod g_i$ are $\tC{p}{m}$.
\item[\thetag{5}] If $f$ is $\tC{p}{m}$, then $t \mapsto f(t^N)$ is $\tC{p}{m}$ for all $N \in \N_{>0}$.
\item[\thetag{6}] If $f$ is $\tC{p}{p}$, then $t \mapsto f(\on{sgn}(t)|t|^{\frac1N})$ is $\tC{p}{p}$ for all $N \in \N_{>0}$.
\item[\thetag{7}] If $f$ is $\tC{p}{p}$, then $t \mapsto \on{sgn}(t)^{l}|t|^\ep f(t)$ is $\tC{p}{p}$ for all $0<\ep<1$ and $l=1,2$.
\end{enumerate}
\end{proposition}

\begin{demo}{Proof}
\thetag{1} and \thetag{2} follow from the definition.

\thetag{3} Cf.\ \cite[Thm.\ 9]{Reichard80}. Clearly $g$ and $F \o g$ are $C^{m-p}$ near $0$ and $C^m$ off $0$. 
By Fa\`a di Bruno's formula \cite{FaadiBruno1855}, for $1 \le k \le p$ and $t\ne0$,
\begin{align} \label{eq:FdB}
  \begin{split}
    \frac{t^k (F \o g)^{(m-p+k)}(t)}{(m-p+k)!} &= \sum_{l\ge 1}  \sum_{\al \in A}
    \frac{t^{k-|\be|}}{l!}  d^l F(g(t)) 
    \Big( 
    \frac{t^{\be_1} g^{(\al_1)}(t)}{\al_1!},\dots,
    \frac{t^{\be_l} g^{(\al_l)}(t)}{\al_l!}\Big) \\
    A &:= \{\al\in \N_{>0}^l : \al_1+\dots+\al_l =m-p+k\} \\
    \be_i &:= \max\{\al_i - m + p,0\}, \quad |\be| = \be_1+\dots+\be_l \le k,
  \end{split}
\end{align}
whose limit as $t\to 0$ exists as a finite number by assumption.

\thetag{4} is a consequence of \thetag{3}.

\thetag{5} As in \eqref{eq:FdB}, for $1 \le k \le p$, $t\ne0$, and $C_{N,\al} := \prod_{i=1}^l \binom{N}{\al_i}$, 
\begin{align*}
\frac{t^k \p_t^{m-p+k}(f(t^N))}{(m-p+k)!} 
&= \sum_{j = -m+p+1}^k  \sum_{_{\substack{\al \in A\\ \al_i\le N\\ l=m-p+j}}}
C_{N,\al} \frac{t^{(N-1)(m-p) +j N}f^{(m-p+j)}(t^N)}{(m-p+j)!}    
\end{align*}
whose limit as $t\to 0$ exists as a finite number by assumption.

\thetag{6} As in \eqref{eq:FdB}, for $1 \le k \le p=m$, $t\ne0$, $N>1$, and $D_{N,\al} := \prod_{i=1}^l \binom{\frac{1}{N}}{\al_i}$,
\begin{align*}
\frac{t^k \p_t^{k}(f(\on{sgn}(t)|t|^{\frac1N}))}{k!} 
&= \sum_{l = 1}^k  \sum_{\al \in A}
D_{N,\al} \frac{\on{sgn(t)}^l |t|^{\frac{l}{N}} f^{(l)}(\on{sgn}(t)|t|^{\frac1N})}{l!}    
\end{align*}
whose limit as $t\to 0$ exists as a finite number by assumption.

\thetag{7} By the Leibniz rule, we have, for $0 \le k \le p$ and $t \ne 0$, 
\begin{align*}
  t^k \p_t^k (\on{sgn}(t)^{l}|t|^\ep f(t)) 
  = \on{sgn}(t)^{l} |t|^{\ep} \sum_{j=0}^k \binom{k}{j} \binom{\ep}{k-j} (k-j)!\, t^j f^{(j)}(t),
\end{align*}
which converges to $0$ as $t \to 0$ by assumption.
The proof is complete.
\qed\end{demo}

\section{A rooted tree associated with \texorpdfstring{$P$}{P}}

\subsection{Derived polynomials} \label{ssec:derP}
Let $P$ be a germ at $0 \in \R$ of a $C^p$-curve of polynomials \eqref{eq:P}, 
where $p \in \N \cup\{\infty\}$. We shall repeatedly use several derived polynomials:
\begin{basedescript}{\desclabelstyle{\pushlabel}\desclabelwidth{2cm}}
  \item[$\quad P=\prod P_i$] Lemma~\ref{split} provides a factorization
    $P=P_1 \cdots P_l$ into germs $P_i$ of $C^p$-curves of polynomials such that the roots of each $P_i(0)$ 
    coincide, but distinct $P_i$ have distinct roots. If $l>1$ we say that \emph{$P$ splits}.
  \item[$\quad \overline P$] $\overline P(z) := P(z+\frac{1}{n}a_1(P))$ is a germ of a $C^p$-curve of polynomials satisfying 
    $a_1(\overline P)=0$.
  \item[$\quad P_{(r)}$] If $a(P)$ is $r$-divisible we define $P_{(r)}$ by setting $a(P_{(r)}) := a(P)_{(r)}$. 
    Then $P_{(r)}$ is $C^{p-nr}$ if $p \ge nr$. 
    If $a(P)$ is strictly $r$-divisible (for instance if $p>nr$) and $a_1(P) = 0$, then 
    $P_{(r)}$ splits. If $\mu_i(t)$ is a choice of the roots for $P_{(r)}(t)$, 
    then $\la_i(t)=t^r\mu_i(t)$ represent the roots of $P(t)$.
  \item[$\quad P^{\pm,N}$] For $N \in \N_{>0}$ we set $P^{\pm,N}(t):=P(\pm t^N)$.
\end{basedescript}

\begin{algorithm}[Local factorization of a curve of polynomials] 
\label{alg}
Let $P$ be a germ at $0 \in \R$ of a $C^p$-curve of polynomials as in \eqref{eq:P}, 
where $p \in \N \cup\{\infty\}$. 
The algorithm will associate with $P$ a rooted tree $T=T(P)=T(P,0)$ whose vertices are labeled by pairs 
$(d,q) \in \N \times \Q$. 
Each vertex in $T$ corresponds uniquely to a (intermediate) factor in the factorization produced by the following steps.
At the beginning $T$ consists just of its root which is labeled $(\deg(P)=n, q)$, 
where $q$ is $0$ if $P$ splits or otherwise $q$ will be determined in \thetag{II}.  

Abusing notation we shall denote by $P$ also an intermediate factor produced in the course of the algorithm.
Let $V$ denote the vertex associated with $P$ and $(d(V)=\deg(P),q(V))$ its label.
\begin{enumerate}[(IIa)]
	\item[(I)] If $P$ splits, $P=\prod P_i$,
		add one vertex $V_i$ for each $P_i$, set $d(V_i) := \deg(P_i)$, and join $V_i$ to $V$ by an edge.
		Feed each $P_i$ into \thetag{II}.
	\item[(II)]	If $P$ does not split, i.e., all roots of $P(0)$ coincide, replace $P$ by $\overline P$. So we can assume that 
    $a_1(P)=0$ and thus all roots of $P(0)$ are equal to $0$. Set 
    \begin{equation} \label{eq:m}
      m = m(P,0):=\min_k \frac{m_0(a_k)}{k}  \in \Q_{\ge 0} \cup \{\infty\}.
    \end{equation}
    If $p< \lceil nm \rceil$ set $q(V) := 0$ and $V$ becomes a leaf of $T$. Otherwise proceed as follows.     
	\item[(IIa)] If $0 \le m < \infty$
	    then $m_0(a_k) \ge km$ for all $k$ and we have equality for some $k$.
      If $t^{-k m}a_{k}(t)|_{t=0}=0$ for all $k$ with $m_0(a_k) = km$,  
      set $q(V):=0$ and $V$ becomes a leaf of $T$. 
      So assume the contrary. Let
      \begin{equation} \label{eq:rN}
        \text{$N\in \N_{>0}$ be minimal so that $r:= Nm \in \N$}.
      \end{equation}
      	  
    If $N=1$, then $a(P)$ is strictly $r$-divisible. 
    Set $q(V) := r$, associate $V$ with $P_{(r)}$, 
    and feed $P_{(r)}$ into \thetag{I}. 
    
		If $N>1$ (equivalently $\frac{r}{N} \not\in \Z$), set $q(V) := \frac{r}{N}$ and  
		$V$ becomes a leaf of $T$. 
		For later use we associate with this leaf the new rooted tree $T(P^{\pm,N},0)$.
    Note that $a(P^{\pm,N})$ is strictly $r$-divisible. 
	\item[(IIb)] If $a(P)=0$ identically,  
		set $q(V) := 0$ and add $n$ vertices, each labeled $(1,0)$ and each joined to $V$ by an edge.
	\item[(IIc)] If $m=\infty$ and $a(P) \ne 0$, 
    set $q(V) := 0$ and add $n$ vertices, each labeled $(1,0)$ and each joined to $V$ by an edge. 
    In this case any continuous root $\la$ of $P$ satisfies $m_0(\la)=\infty$.
    Indeed, $a(P)$ is $r$-divisible for every $r \in \N$, and, for $t\ne 0$, $\mu(t) = t^{-r} \la(t)$ is a root of $P_{(r)}(t)$, 
    thus bounded, and so $m_0(\la) \ge r-1$.
\end{enumerate} 
\end{algorithm}

\subsection{Good, bad, and fair rooted trees \texorpdfstring{$T(P,t)$}{T(P,t)}} \label{gbf}
Modifying the algorithm in the obvious way we can associate with each $C^p$-curve of polynomials 
$P(t)$, $t \in I$, as in \eqref{eq:P} and each $t \in I$ a rooted tree 
$T(P,t)$ whose vertices are labeled by pairs $(d,q) \in \N \times \Q$.

We say that $T(P,t)$ is \emph{admissible} if the label $(d,q)$ of every leaf of $T(P,t)$ satisfies
\begin{itemize}
  \item $q = 0 \Rightarrow d=1$.
\end{itemize}
That is, in Step~\thetag{II} of Algorithm~\ref{alg}, $m$ satisfies $0 < m \le \frac{p}{n}$ and $a(P)$, resp.\ $a(P^{\pm,N})$, 
is strictly $r$-divisible.
We say that $T(P,t)$ is \emph{good} if every leaf of $T$ has the label $(1,0)$ (so a good tree is admissible). 
An admissible tree $T(P,t)$ that is not good is called \emph{bad}. Leaves with label other than $(1,0)$ are called bad.
An admissible tree $T(P,t)$ is called \emph{fair} if 
for each bad leaf $V$ in $T(P,t)$ the rooted tree $T(P_V^{\pm,N},t)$, where $P_V$ is the polynomial corresponding to $V$
(see Step~\thetag{IIa}), is good and has height~$1$.

For $P$ hyperbolic and $p$ sufficiently large, $T(P,t)$ is good for all $t$, see e.g.\  \cite{RainerOmin}.

\subsection{Associated integers \texorpdfstring{$\Ga(P)$}{Ga(P)}, \texorpdfstring{$\ga(P)$}{ga(P)}, and 
\texorpdfstring{$\hat \Ga(P)$}{hatGaP}} \label{Gaga}
Let $P(t)$, $t \in I$, be a $C^p$-curve of polynomials as in \eqref{eq:P}.
For $t \in I$ and each vertex $V$ in $T(P,t)$ we recursively define two integers $\Ga_{t}(V)$ and $\ga_{t}(V)$ by setting
\begin{align*}
  \Ga_{t}(V) &:= \lceil d(V) \cdot q(V) \rceil \\
    \ga_{t}(V) &:= \lfloor q(V) \rfloor 
\end{align*}
if $V$ is a leaf, and otherwise 
\begin{gather*} 
    \Ga_{t}(V) := \lceil d(V) \cdot q(V) \rceil + \max_W \Ga_{t}(W), \\
    \ga_{t}(V) := \lfloor q(V) \rfloor +
    \left\{
    \begin{array}{ll}
    \Ga_{t}(V) - \max_W \big(\Ga_{t}(W)-\ga_{t}(W)\big) & \quad\text{if } T(P,t) \text{ is good} \\
    \min_W \ga_{t}(W) & \quad\text{if } T(P,t) \text{ is bad}
    \end{array}
    \right.,
\end{gather*}
where $W$ ranges over all successors of $V$, and we set
\begin{align} \label{eq:alt}
  \begin{split}
    \Ga_{t}(P) &:= \Ga_{t}(\text{root of }T(P,t)), \\
    \ga_{t}(P) &:= \ga_{t}(\text{root of }T(P,t)).
  \end{split}
\end{align}
Furthermore, we set 
\begin{gather}
\Ga(P) := \sup_{t \in I} \Ga_{t}(P), \label{Ga_sup}\\
\ga(P) := 
\min\Big\{\inf_{\substack{t \in I:\\ T(P,t) \text{ good}}} \big(\Ga(P) -\Ga_{t}(P) + \ga_{t}(P)\big),
\inf_{\substack{t \in I:\\ T(P,t) \text{ bad}}} \ga_{t}(P)\Big\} \label{ga_inf}.
\end{gather}
Then $\Ga(P),\ga(P) \in \N \cup \{\infty\}$ and $\ga(P) \le \Ga(P)$. 

\medskip

For later reference we also present an equivalent definition of $\Ga_{t}(P)$ and $\ga_{t}(P)$ under the assumption 
that $T(P,t)$ is admissible:
\begin{itemize}
	\item
    If $P = \prod P_i$ splits near $t$,  
		\begin{gather}
		\Ga_{t}(P) := \max_i \Ga_{t}(P_i), \label{Ga_split}\\
		\ga_{t}(P) := 
		\left\{
		\begin{array}{ll}
		\Ga_{t}(P) - \max_i \big(\Ga_{t}(P_i)-\ga_{t}(P_i)\big) & \quad\text{if } T(P,t) \text{ is good} \\
		\min_i \ga_{t}(P_i) & \quad\text{if } T(P,t) \text{ is bad}
		\end{array}
		\right.
		. \label{ga_split}
		\end{gather}
  \item 
    If $P$ does not split near $t$, 
    \begin{align}
    \Ga_{t}(P) &:= \Ga_{t}(\overline P),  \label{Ga_bar}\\
    \ga_{t}(P) &:= \ga_{t}(\overline P). 
    \end{align} 
    Let $m$ be as defined in \eqref{eq:m}; as $T(P,t)$ is admissible $m\ne 0$.
    \begin{itemize}
      \item If $m <\infty$ consider the integers $N$ and $r$ from \eqref{eq:rN}.
    \begin{itemize}
      \item If $N=1$, 
    \begin{align}
    \Ga_{t}(P) &:= \Ga_{t}(P_{(r)}) + \deg(P) r, \label{Ga_(r)}\\
    \ga_{t}(P) &:= \ga_{t}(P_{(r)})+ r. \label{ga_(r)}
    \end{align} 
      \item If $N>1$, 
    \begin{align}
    \Ga_{t}(P) &:= \Big\lceil \frac{\deg(P) r}{N} \Big\rceil, \label{GaN>1}\\
    \ga_{t}(P) &:= \Big\lfloor \frac{r}{N} \Big\rfloor. \label{gaN>1}
    \end{align}
    \end{itemize}
    \item If $m =\infty$,  
    \begin{align}
      \Ga_{t}(P) := 0, \label{Ga_infty} \\ 
      \ga_{t}(P) := 0.
    \end{align}      
    \end{itemize}
\end{itemize}

Still assuming that $T(P,t)$ is admissible, we define
$\hat \Ga_{t}(P)$, resp.$\ \hat \Ga(P)$, by replacing $\Ga$ by $\hat \Ga$ in 
\eqref{Ga_sup}, \eqref{Ga_split}, \eqref{Ga_bar}, \eqref{Ga_(r)}, \eqref{Ga_infty} and
by substituting 
\begin{equation} \label{hatGaN>1}
  \hat \Ga_{t}(P):= \hat \Ga_{t}(P^{\pm,N}).
\end{equation}
for \eqref{GaN>1}. Clearly $\Ga_t(P) \le \hat \Ga_t(P)$ and $\Ga(P) \le \hat \Ga(P)$. 

\begin{remark} \label{rem:admissible}
We have:
  \begin{enumerate}
    \item If $T(P,t)$ is admissible, then $P$ is $C^{\Ga_{t}(P)}$ near $t$.
    \item If $P$ is $C^{\Ga_{t}(P)+1}$ near $t$, then $T(P,t)$ is admissible. 
  \end{enumerate}
\end{remark}

\subsection{The set \texorpdfstring{$E^{(\infty)}(P)$}{}}
Let $P(t)$, $t \in I$, be a $C^0$-curve of polynomials as in \eqref{eq:P}. 
We denote by $E^{(\infty)}(P)$ the set of all $t \in I$ which satisfy the following condition: 
\begin{itemize}
\item If $s=s(t,P)$ is the maximal integer so that the germ at $t$ of $\tilde \De_s \o P$ is not $0$,
then $m_t(\tilde \De_s \o P) = \infty$.
\end{itemize}
Let $\La_t$ denote the set of germs at $t$ of $C^0$-parameterizations $\la = (\la_1,\ldots,\la_n)$ of the roots of $P$ which is 
non-empty by Lemma~\ref{controots}. 
If $t \not\in E^{(\infty)}(P)$ then $m_t(\la_i-\la_j) = \infty$ implies $\la_i=\la_j$ for all $\la \in \La_t$, by \eqref{eq:del}, 
and Step~\thetag{IIc} in Algorithm~\ref{alg} is empty.

\section{Differentiable roots of complex polynomials}

In this section we prove the main result, Theorem~\ref{main}. First we state two lemmas.

\begin{lemma}[{\cite[2.6]{RainerN}}] \label{unique}
  Let $P(t)$, $t \in I$, be a $C^0$-curve of polynomials as in \eqref{eq:P} such that $E^{(\infty)}(P)= \emptyset$.
  Two $C^\infty$-parameterizations of the roots of $P$ differ by a constant permutation.  
\end{lemma}

The following lemma is needed for polynomials with fair rooted tree.

\begin{lemma} \label{bad}
Let $s, N \in \N_{>0}$ be such that $q:=\frac{s}{N} \not\in \Z$ and let $r \in \N$. 
Let $P$ be a germ at $0$ of a $\tC{nr}{\lceil n q \rceil + nr}$-curve of polynomials as in \eqref{eq:P} such that
\begin{itemize}
  \item $a_1(P)=0$
  \item $m_0(a_k(P)) \ge kq$ for all $k$ (thus $a(P^{\pm,N})$ is $s$-divisible)
  \item all roots of $(P^{\pm,N})_{(s)}(0)$ are distinct.
\end{itemize}
Let $\mu^{\pm,N}_j$ be continuous function germs representing the roots of $(P^{\pm,N})_{(s)}$ and set $\la^{\pm,N}_j(t):= t^s \mu^{\pm,N}_j(t)$.
Then each $\La^{\pm}_j(t) := \la^{\pm,N}_j(\on{sgn}(t) |t|^{\frac1N})$ is $\tC{r}{\lfloor q \rfloor+r}$.
\end{lemma}

\begin{demo}{Proof}
Write $a_k=a_k(P)$.
Clearly $m_0(a_k) \ge \lceil k q \rceil$ and thus
$a_k(t) = t^{\lceil k q \rceil} b_k(t)$ for $\tC{\lceil n q \rceil +nr}{\lceil n q \rceil + nr}$-functions $b_k$, 
by \ref{prop}\thetag{2}.
By \ref{prop}\thetag{5}, $t \mapsto b_k(\pm t^N)$ is $\tC{\lceil n q \rceil +nr}{\lceil n q \rceil + nr}$ as well.
It follows that $a((P^{\pm,N})_{(s)})$ is $\tC{\lceil n q \rceil +nr}{\lceil n q \rceil + nr}$, since
\[
a_k((P^{\pm,N})_{(s)})(t) = (\pm 1)^{\lceil k q \rceil} t^{N \lceil k q \rceil-ks} b_k(\pm t^N), \quad N \lceil k q \rceil-ks \in \N.  
\]
By Lemma~\ref{split} there exist real analytic functions $\Ph_{j}$, $1 \le j \le n$, defined in a neighborhood of 
$a((P^{\pm,N})_{(s)})(0) \in \C^{n}$ such that $\mu^{\pm,N}_{j} = \Ph_{j} \o a((P^{\pm,N})_{(s)})$.
Hence each $t \mapsto \mu^{\pm,N}_{j}(\on{sgn}(t)|t|^{\frac1N})$ is $\tC{\lceil n q \rceil +nr}{\lceil n q \rceil + nr}$, 
by \ref{prop}\thetag{3} 
and \ref{prop}\thetag{6}.
Now
\begin{equation} \label{eq:Lambda}
\La^{\pm}_j(t) 
= \on{sgn}(t)^{s+\lfloor q \rfloor} |t|^{\{q\}} t^{\lfloor q \rfloor} \mu^{\pm,N}_{j}(\on{sgn}(t)|t|^{\frac1N})
\end{equation}
is $\tC{r}{\lfloor q \rfloor + r}$, by \ref{prop}\thetag{1} and \ref{prop}\thetag{7}, since $\{q\}\ne 0$. 
\qed\end{demo}

Now we are ready to prove the main theorem.

\begin{theorem} \label{main}
Let $P(t)$, $t \in I$, be a $C^{p+\Ga(P)}$-curve, $p \in \N \cup \{\infty\}$, of polynomials as in \eqref{eq:P} 
such that $E^{(\infty)}(P)=\emptyset$.
Then we have:
\begin{enumerate}
\item[\thetag{1}] If $T(P,t)$ is good for all $t$, 
then $P$ is $C^{p+\ga(P)}$-solvable.
\item[\thetag{2}] If $T(P,t)$ is fair for all $t$,
then $P$ is $C^{\ga(P)}$-solvable.
\item[\thetag{3}] If $T(P,t)$ is bad for some $t$ and $P$ is $C^{1+\hat \Ga(P)}$,
then $P$ is $C^{\ga(P)}$-solvable.
\end{enumerate}
\end{theorem}

Note that $T(P,t)$ is admissible for all $t$; in \thetag{1} and \thetag{2} by definition and in \thetag{3} 
by Remark~\ref{rem:admissible}.

\begin{demo}{Proof}
It suffices to prove the local assertions in the following claim. 
This is a consequence of Lemma~\ref{glue} and Lemma~\ref{unique}.
Let $t$ be fixed.

\begin{claim*}[Local assertions]
Assume $P$ is $C^{p+\Ga(P)}$.
\begin{enumerate} 
\item[\thetag{$1'$}] If $T(P,t)$ is good, then $P$ is $C^{p+\ga(P)}$-solvable near $t$. 
\item[\thetag{$2'$}] If $T(P,t)$ is fair, then $P$ is $C^{\ga(P)}$-solvable near $t$.
\item[\thetag{$3'$}] If $T(P,t)$ is bad and $P$ is $C^{1+\hat \Ga(P)}$, then $P$ is 
$C^{\ga(P)}$-solvable near $t$.  
\end{enumerate}
\end{claim*}

In fact it suffices to show the following.

\begin{claim*}
Assume that $P$ is $C^{p+\Ga_{t}(P)}$.
\begin{enumerate}
\item[\thetag{$1''$}] If $T(P,t)$ is good, 
then $P$ is $C^{p+\ga_{t}(P)}$-solvable near $t$. 
\item[\thetag{$2''$}] If $T(P,t)$ is fair, then $P$ is $C^{\ga_{t}(P)}$-solvable near $t$.
\item[\thetag{$3''$}] If $T(P,t)$ is bad and $P$ is $C^{1+\hat \Ga_{t}(P)}$, 
then $P$ is $C^{\ga_{t}(P)}$-solvable near $t$.  
\end{enumerate}
\end{claim*}

Indeed, if $P$ is $C^{p+\Ga(P)}$ and we assume \thetag{$1''$}, then $P$ is  
$C^{p+\Ga(P)-\Ga_{t}(P) + \ga_{t}(P)}$-solvable near $t$, and, by \eqref{ga_inf},  
$p+\Ga(P)-\Ga_{t}(P) + \ga_{t}(P) \ge p + \ga(P)$, thus \thetag{$1'$}. The implications
\thetag{$2''$} $\Rightarrow$ \thetag{$2'$} and \thetag{$3''$} $\Rightarrow$ \thetag{$3'$} follow from 
$\ga_{t}(P) \ge \ga(P)$, by \eqref{ga_inf}.

\medskip

Without loss of generality assume that $0 \in I$ and $t=0$.
We obtain \thetag{$1''$}, \thetag{$2''$}, and \thetag{$3''$} if we set $r=0$ and $P_{[0]}=P$ in the following claim.

\begin{claim*}
Let $P_{[r]}$ be any $\tC{nr}{p+\Ga_0(P_{[r]}) + nr}$-curve, $r \in \N$, of polynomials as in \eqref{eq:P} such that 
$0 \not\in E^{(\infty)}(P_{[r]})$.    
\begin{enumerate}
\item[\thetag{$1'''$}] If $T(P_{[r]},0)$ is good, 
then $P_{[r]}$ is $\tC{r}{p+\ga_0(P_{[r]})+r}$-solvable. 
\item[\thetag{$2'''$}] If $T(P_{[r]},0)$ is fair,
then $P_{[r]}$ is $\tC{r}{\ga_0(P_{[r]})+r}$-solvable. 
\item[\thetag{$3'''$}] If $T(P_{[r]},0)$ is bad and $P_{[r]}$ is $\tC{nr}{1+\hat \Ga_0(P_{[r]}) + nr}$,
then $P_{[r]}$ is $\tC{r}{\ga_0(P_{[r]})+r}$-solvable. 
\end{enumerate}
\end{claim*}

Here $P_{[r]}$ should not be confused with $P_{(r)}$; 
however, if $P \leadsto P_{(r)}$ as in Step~\thetag{IIa} of Algorithm~\ref{alg}, 
then $P_{(r)}$ satisfies the assumptions of the claim. 

We first prove \thetag{$1'''$} and \thetag{$2'''$} by induction on the degree $n$; afterwards we indicate the modifications 
which provide a proof of \thetag{$3'''$}. 
We follow Algorithm~\ref{alg}.

\smallskip

\begin{demo}{Proof of ($1'''$) and ($2'''$)}
	If all roots of $P_{[r]}(0)$ coincide, we may assume that $a_1(P_{[r]})=0$, 
	by replacing $P_{[r]}$ by $\overline{P_{[r]}}$ which preserves the class $\tC{nr}{p+\Ga_0(P_{[r]}) + nr}$,   
	by \ref{prop}\thetag{4}. 
	Then all roots of $P_{[r]}(0)$ are equal to $0$,
	and $a(P_{[r]})(0)=0$.

	If $a(P_{[r]})=0$ identically, we are done.
  Otherwise consider 
  $m = m(P_{[r]},0) \in \Q_{>0}$ as in \eqref{eq:m} 
  and let $N \in \N_{>0}$ be minimal such that $s:=Nm \in \N_{>0}$; note that $0<m<\infty$ since $T(P_{[r]},0)$ is admissible
  and $0 \not\in E^{(\infty)}(P_{[r]})$.  
  We treat the cases $N=1$ and $N>1$ separately.

  \smallskip

	\textbf{Case $N=1$.} 
  Admissibility of $T(P_{[r]},0)$ implies that $a(P_{[r]})$ is strictly $s$-divisible.
  The polynomial $P_{(r+s)} := (P_{[r]})_{(s)}$ is $\tC{n(r+s)}{p+\Ga_0(P_{[r]}) + nr}$, 
  by \ref{prop}\thetag{2}, and it splits.  
  Note that $T(P_{[r]},0) = T(P_{(r+s)},0)$, $\Ga_0(P_{[r]})=\Ga_0(P_{(r+s)}) + ns$, by \eqref{Ga_(r)}, 
  $\ga_0(P_{[r]})=\ga_0(P_{(r+s)}) + s$, by \eqref{ga_(r)}, $0 \not\in E^{(\infty)}(P_{(r+s)})$, 
  and the roots of $P_{[r]}$ differ from those of $P_{(r+s)}$ by multiplication with $t^s$. 
	Replace $P_{[r]}$ by $P_{(r+s)}$. 

	So after these initial steps we may assume that $a_1(P_{[r]})=0$ and that $P_{[r]}$ splits:
	$P_{[r]} = \prod P_{[r],j}$, where $P_{[r],j} := (P_{[r]})_j$ and 
	$a(P_{[r],j}) = \Ph_j \o a(P_{[r]})$ for real analytic mappings $\Ph_j$ defined near $a(P_{[r]})(0) \in \C^{n}$, 
	by Lemma~\ref{split}.
	Each $P_{[r],j}$ forms a $\tC{nr}{p+\Ga_0(P_{[r]}) + nr}$-curve of polynomials, by \ref{prop}\thetag{3}.
	Setting 
  \[
    p_{[r],j}:=p+\Ga_0(P_{[r]}) -\Ga_0(P_{[r],j}),
  \] 
  we may conclude that
	$P_{[r],j}$ is $\tC{n_j r}{p_{[r],j}+\Ga_0(P_{[r],j}) + n_j r}$, where $n_j=\deg(P_{[r],j})$, by \ref{prop}\thetag{1}. 
	If $T(P_{[r]},0)$ is good (resp.\ fair), then each $T(P_{[r],j},0)$ is good (resp.\ fair).
	If $T(P_{[r]},0)$ is bad, then some $T(P_{[r],j},0)$ is bad. 
	Clearly 
	$0 \not\in E^{(\infty)}(P_{[r],j})$ for all $j$. Hence each $P_{[r],j}$ satisfies 
	the assumptions of the claim and thus is 
	\begin{enumerate}
	\item[\thetag{$1'''$}] $\tC{r}{p_{[r],j}+\ga_{0}(P_{[r],j}) + r}$-solvable,
	\item[\thetag{$2'''$}] $\tC{r}{\ga_{0}(P_{[r],j}) + r}$-solvable,
	\end{enumerate}
  by induction.
	By \eqref{Ga_split} and \eqref{ga_split}, we have $p_{[r],j}+\ga_{0}(P_{[r],j}) \ge p+\ga_{0}(P_{[r]})$ 
  (resp.\ $\ga_{0}(P_{[r],j}) \ge \ga_{0}(P_{[r]})$) 
	if $T(P_{[r]},0)$ is good (resp.\ bad).
	Hence we have shown \thetag{$1'''$} and \thetag{$2'''$} in case $N=1$.

  \smallskip

	\textbf{Case $N>1$.} Here $T(P_{[r]},0)$ is trivial, the label of the only vertex is  
  $(n,q)$ where $q=\frac{s}{N}$. 
  Admissibility of $T(P_{[r]},0)$ implies that $a((P_{[r]})^{\pm,N})$ is strictly $s$-divisible.
  By Lemma~\ref{controots}, $P_{(r+s)}^{\pm,N} := ((P_{[r]})^{\pm,N})_{(s)}$ 
	admits a continuous parameterization $\mu^{\pm,N}_j$ of its roots satisfying $\mu^{+,N}_j(0)=\mu^{-,N}_j(0)$, 
	and $\la^{\pm,N}_j(t) := t^s \mu^{\pm,N}_j(t)$ represent the roots of $(P_{[r]})^{\pm,N}$.
	As $T(P_{[r]},0)$ is fair, $\mu^{\pm,N}_i(0) \ne \mu^{\pm,N}_j(0)$ for $i \ne j$, and thus each 
	$\La^{\pm}_j(t) := \la^{\pm,N}_j(\on{sgn}(t) |t|^{\frac1N})$ is $\tC{r}{\lfloor q \rfloor+r}$, by Lemma~\ref{bad}.
	So, if $N$ is odd, then $\la_j:= \La^{+}_j$ forms a $\tC{r}{\lfloor q \rfloor+r}$-parameterization of the roots of $P_{[r]}$.
	If $N$ is even, set 
	\begin{equation*} \label{eq:ref}
	\la_j(t):= 
	\begin{cases}
	\La^{+}_j(t) & \quad\text{if } t\ge 0\\
	\La^{-}_j(t) & \quad\text{if } t< 0
	\end{cases}
	.
	\end{equation*}
	The functions $\la_j$ parameterize the roots of $P_{[r]}$ and are $\tC{r}{\lfloor q \rfloor+r}$. 
  Indeed, using the formula \eqref{eq:Lambda} for $\La^{\pm}_j$ and 
  applying the Leibniz rule, we find  
  $\p_t^m(t^r \La^{\pm}_j(t))|_{t=0} = 0$ for all $0 \le m \le \lfloor q \rfloor + r$, 
	since $\{q\} \ne 0$.  
	By \eqref{gaN>1}, we are done. The proof of ($1'''$) and ($2'''$) is complete.	
\end{demo}

\begin{demo}{Proof of ($3'''$)}
	The first part of the proof of \thetag{$2'''$} works for \thetag{$3'''$} as well; just exchange $\hat \Ga_0$ for $\Ga_0$ and set $p=1$.
	The case $N>1$ must be modified:

	If $T(P_{[r]},0)$ is not fair, 
	the roots $\mu^{\pm,N}_j$ of $P_{(r+s)}^{\pm,N}$ are not necessarily pairwise distinct at $0$ and we cannot apply 
	Lemma~\ref{bad}. 
	However we may proceed as follows.
	By \ref{prop}\thetag{5}, $(P_{[r]})^{\pm,N}$ is $\tC{nr}{1+\hat \Ga_0(P_{[r]}) + nr}$, 
	and $0 \not\in E^{(\infty)}((P_{[r]})^{\pm,N})$.
	Since $\hat \Ga_0(P_{[r]}) = \hat \Ga_0((P_{[r]})^{\pm,N})$, by \eqref{hatGaN>1}, 
	and since $a((P_{[r]})^{\pm,N})$ is strictly $s$-divisible, 
  we may conclude by induction and the first part of the proof 
  that the roots of $(P_{[r]})^{\pm,N}$ can be parameterized by 
	$\tC{r}{\ga_0((P_{[r]})^{\pm,N})+r}$-functions $\la_j^{\pm,N}$.   
	Since $\la^{\pm,N}_j(t) = t^s \mu^{\pm,N}_j(t)$ and by \eqref{ga_(r)} and \ref{prop}\thetag{1}, 
	we find that $\mu^{\pm,N}_j$ is $\tC{r+s}{r+s}$. 
	By \ref{prop}\thetag{6}, \ref{prop}\thetag{1}, and \ref{prop}\thetag{7}, we may conclude that 
	$\La^{\pm}_j(t) = \la^{\pm,N}_j(\on{sgn}(t) |t|^{\frac1N})$ is $\tC{r}{\lfloor q \rfloor+r}$
  (cf.\ the end of the proof of Lemma~\ref{bad}). 
	We finish the proof as for \thetag{$2'''$}.	
\end{demo}

The proof is complete.
\qed\end{demo}

\begin{remarks}
  \thetag{1}
  If $T(P,t)$ is admissible, then so is the tree associated with any of the derived polynomials of $P$, except possibly for $P^{\pm,N}$ 
  (cf.\ \ref{ssec:derP}). 
  The condition that $P$ is $C^{1+\hat \Ga(P)}$ in \ref{main}\thetag{3} guarantees admissibility of the associated trees of $P$ 
  at any $t$
  and of all iterated derivations including $(\cdot)^{\pm,N}$ which appear in the course of Algorithm~\ref{alg}, see 
  Remark~\ref{rem:admissible}. If admissibility 
  of all these trees is taken for granted, 
  then the condition that $P$ is $C^{\hat \Ga(P)}$ is sufficient for the conclusion of \ref{main}\thetag{3}.  

  \thetag{2} A priori 
$\Ga(P)$, $\hat \Ga(P)$, and $\ga(P)$ might be infinite. Theorem~\ref{main} remains true in this case; 
note that $\ga(P)$ is finite if $T(P,t)$ is bad for some $t$. 
The integers $\Ga_t(P)$, $\hat \Ga_t(P)$, and $\ga_t(P)$ are non-zero only at points $t$, where the multiplicity of the roots of $P$
changes, and these points form a discrete set if $E^{(\infty)}(P)= \emptyset$ and $P$ is $C^\infty$. 
So $\Ga_J(P)$, $\hat \Ga_J(P)$, and $\ga_J(P)$, 
where the supremum, resp.\ infimum, in the definition \eqref{Ga_sup} and \eqref{ga_inf} is taken over a relatively compact subinterval 
$J \subseteq I$, are finite if $E^{(\infty)}(P)= \emptyset$ and $P$ is $C^\infty$.
\end{remarks}

\begin{corollary} \label{cor1}
Let $p \in \N \cup \{\infty\}$ and let $f$ be a germ at $0 \in \R$ of a complex valued $C^{p+m_0(f)}$-function 
with $m_0(f) <\infty$ and so that $t^{-m_0(f)}f(t)|_{t=0} \ne 0$.  
Then:
\begin{enumerate}
\item[\thetag{1}] If $\frac{m_0(f)}{n} \in \Z$, then there is a $C^{p+\frac{m_0(f)}{n}}$-germ $g$ such that $g^n=f$.
\item[\thetag{2}] If $\frac{m_0(f)}{n} \not\in \Z$, then there is a $C^{\lfloor \frac{m_0(f)}{n} \rfloor}$-germ $g$ 
such that $g^n=f$.
\end{enumerate}
\end{corollary}

\begin{demo}{Proof}
Consider $P(t)(z) = z^n-f(t)$. Then $\frac{r}{N}=\frac{m_0(f)}{n}$ and $\Ga_0(P) = m_0(f)$.
Thus, \thetag{1} is a special case of \ref{main}\thetag{$1''$}, and
\thetag{2} of \ref{main}\thetag{$2''$}.
\qed\end{demo}

This result is essentially due to Reichard \cite{Reichard80}.
Note that there are minor differences between Corollary~\ref{cor1} and \cite[Cor.\ 15]{Reichard80}, since 
Reichard considers the $n$th root $f^{\frac{1}{n}}$ of a (non-negative if $n$ is even) function $f$ while we study
functions $g$ such that $g^n=f$.

\section{Differentiable roots of definable complex polynomials} \label{sec:def}

For definable curves of polynomials Theorem~\ref{main} (and its proof) simplifies. 
Admissibility of the associated trees as well as
the assumption $E^{(\infty)}(P) = \emptyset$    
are not necessary. The main reason is the following lemma. 

\begin{lemma}[{\cite[2.5]{RainerOmin}}] \label{deflem}
Let $f$ be a germ at $0 \in \R$ of a definable complex valued function and
let $m \in \N$. If $f$ is $C^m$, 
then $f$ is $\tC{p}{m+p}$ for every $p \in \N$. 
In particular, if $f$ is $C^0$ and $m_0(f)\ge p$, then $f$ is $C^p$. 
\end{lemma}

Let $P(t)$, $t \in I$, be a definable $C^0$-curve of polynomials as in \eqref{eq:P}.
Then any continuous parameterization of the roots of $P$ is definable (cf.\ \cite[7.2]{RainerOmin}).
We slightly modify Algorithm~\ref{alg}:
\begin{enumerate}[(IIa)]
  \item[(II)] If $p < \lceil nm \rceil$ proceed in \thetag{IIa}.
  \item[(IIa)] If $t^{-km}a_k(t)|_{t=0} = 0$ for all $k$ with $m_0(a_k) = km$, 
  consider $N,r$ as in \eqref{eq:rN}, set $q(V):=\frac{r}{N}$, and $V$ becomes a leaf of $T$. 
  In this case $a(P)$, resp.\ $a(P^{\pm,N})$, is $r$-divisible, but not strictly.   
\end{enumerate}
Define $\Ga_{t}(P)$ and $\ga_{t}(P)$ by \eqref{eq:alt}; it does not require admissibility of $T(P,t)$.
Definability implies that the set of $t \in I$ such that $\Ga_{t}(P)>0$ or $\ga_{t}(P)>0$ is finite and hence
$\Ga(P)$ and $\ga(P)$ are integers.

\begin{theorem} \label{maindef}
Let $P(t)$, $t \in I$, be a definable $C^{p+\Ga(P)}$-curve, $p \in \N \cup \{\infty\}$, of polynomials as in \eqref{eq:P}.
Then:
\begin{enumerate}
\item[\thetag{1}] If $T(P,t)$ is good for all $t$, 
then $P$ is $C^{p+\ga(P)}$-solvable. 
\item[\thetag{2}] If $T(P,t)$ is bad for some $t$, 
then $P$ is $C^{\ga(P)}$-solvable. 
\end{enumerate}
\end{theorem}

\begin{demo}{Proof}
Here Step~\thetag{IIc} of Algorithm~\ref{alg} might not be empty. 
Suppose that $0 \in I$ and $m_0(a_k) = \infty$ for all 
$k$, where $a_k = a_k(P)$. Then $m_0(\la_j) = \infty$ for any choice of continuous roots $\la_j$ of $P$.
By Lemma~\ref{deflem}, for each $q$ 
there is a neighborhood $I_q$ of $0$ such that each $\la_j$ is $C^q$ on $I_q$.
Since the coefficients $a_j$ (and hence the functions $\tilde \De_k(P)$) are definable,
the multiplicity of the $\la_j(t)$ is constant for small $t \ne 0$. 
So off $0$ all $\la_j$ have the regularity of $P$, by Lemma~\ref{split}, and hence 
also near $0$.

This observation together with \ref{main}\thetag{$1$} implies the case $p=\infty$, since by Lemma~\ref{unique} 
the local $C^\infty$-parameterizations are unique up to flat contact in the following sense:
  If $\la_j$ and $\mu_j$ are two different $C^\infty$-parameterizations of the roots near $0$, then 
  $\{\la_1,\ldots,\la_n\}/\!{\sim}=\{\mu_1,\ldots,\mu_n\}/\!{\sim}$, where $\la_i \sim \la_j$ if and only if 
  $m_0(\la_i-\la_j)=\infty$.  
 
The rest of the proof is similar to the proof of Theorem~\ref{main} but much simpler via Lemma~\ref{deflem}; 
here we automatically ``gain $q$ derivatives back by multiplying with $t^q$'', in particular, 
we do not need ``strict divisibility''.
It is obvious how to make the necessary modifications.
\qed\end{demo}

For the sake of completeness we state the following result.

\begin{theorem} \label{C^1}
Let $P(t)$, $t \in I$, be a definable $C^n$-curve of polynomials
of degree $n$ as in \eqref{eq:P}.
Then $P$ is $C^1$-solvable if and only if there is a continuous parameterization 
of the roots with order of contact $\ge 1$
(i.e.\ if $\la_i(t)=\la_j(t)$ then $m_{t}(\la_i-\la_j) \ge 1$).
\end{theorem}

\begin{demo}{Proof}
The statement follows from Lemma~\ref{glue}, Lemma~\ref{deflem}, and 
\cite[4.3]{RainerAC}.
\qed\end{demo}

\section{Differentiable eigenvalues and eigenvectors of normal matrices} \label{sec:normal}

We say that $A(t)=(A_{ij}(t))_{1 \le i,j \le n}$, $t \in I$, is a \emph{$C^p$-curve of normal matrices}, 
$p \in \N \cup \{\infty\}$, if all $A_{ij} \in C^p(I,\C)$ and 
$A(t)A^*(t)=A^*(t)A(t)$ for all $t$.
We associate with $A$ its characteristic polynomial  
$P_A(z) := \det (z\I_n-A)$, i.e., $a_k(P_A) = \on{Trace}(\bigwedge^k\! A)$, and set 
$E^{(\infty)}(A):=E^{(\infty)}(P_A)$.

\begin{algorithm}[Local spectral decomposition of a curve of normal matrices] \label{alg:n}
  Let $A=(A_{ij})_{1 \le i,j \le n}$ be a germ at $0 \in \R$ of a  
  $C^p$-curve of normal matrices,  
  where $p \in \N \cup\{\infty\}$. 
  In analogy with Algorithm~\ref{alg}, this algorithm will associate with $A$ a rooted tree $T=T(A)=T(A,0)$ 
  whose vertices are labeled by pairs $(d,q) \in \N \times \Q$. 
  At the beginning $T$ consists just of its root which is labeled $(n, q)$, 
  where $q$ is $0$ if we start in \thetag{I} or otherwise $q$ will be determined in \thetag{II}.  

  Let $V$ denote the vertex associated with $A$ and $(d(V),q(V))$ its label.
  \begin{enumerate}[(I)]
  	\item[(I)] If $A(0)$ has distinct eigenvalues $\nu_1,\ldots,\nu_l$ with respective multiplicities 
  		$n_1,\ldots,n_l$, 
		choose pairwise disjoint simple closed $C^1$-curves $\ga_i$ in the resolvent set of $A(0)$ 
		such that $\ga_i$ encloses only $\nu_i$ among
		all eigenvalues of $A(0)$.
		By continuity, no eigenvalue 
		of $A(t)$ lies on $\bigcup_{i} \ga_i$, for $t$ near $0$.
		Now,
		\begin{equation*}
		t\mapsto -\frac1{2\pi i}\int_{\ga_i} (A(t)-z)^{-1}\;dz =: \Pi_i(t)
		\end{equation*}
		is a $C^p$-curve of projections with constant rank onto the direct sum of all 
		eigenspaces corresponding to eigenvalues of $A(t)$ in the interior of $\ga_i$.
		The family of $n_i$-dimensional complex vector spaces 
		$t\mapsto \Pi_i(t)(\C^n)\subseteq \C^n$ forms a $C^p$ Hermitian vector subbundle 
		of the trivial bundle $\R \times \C^n \to \R$:
		For given $t$, choose $v_1,\dots v_{n_i}\in \C^n$ such that the $\Pi_i(t)(v_i)$ are linearly independent 
		and thus span $\Pi(t)(\C^n)$. 
		This remains true locally in $t$. We use the Gram Schmidt orthonormalization 
		procedure (which is $C^\om$) 
		for the $\Pi_i(t)(v_i)$ to obtain a local orthonormal $C^p$-frame of the bundle. 
		Now $A(t)$ maps $\Pi_i(t)(\C^n)$ to itself and in a local $C^p$-frame it is given by a normal $n_i \times n_i$ matrix $A_i$
		parameterized in a $C^p$-way by $t$. 
    Add one vertex $V_i$ for each $A_i$, join $V_i$ to $V$ by an edge, and set $d(V_i) := n_i$.
		Feed each $A_i$ into \thetag{II}.
	\item[(II)] If all eigenvalues of $A(0)$ coincide, replace $A$ by $A - \frac{1}{n} \on{Trace}(A) \I_n$.
        Then all eigenvalues of $A(0)$ are equal to $0$, and $A(0)=0$.	
        Set 
        \begin{equation} \label{eq:r}
          r:= \min_{i,j} m_0(A_{ij}).
        \end{equation}
        If $r=\infty$ set $q(V) := 0$ and add $n$ vertices, each joined to $V$ by an edge and each labeled $(1,0)$. 
        If $A=0$ all eigenvalues of $A$ are identically $0$ and the eigenvectors can be chosen constant.  
        If $A \ne 0$ (which may occur if $0 \in E^{(\infty)}(A)$) then any continuous eigenvalue $\la$ of $A$ 
        satisfies $m_0(\la) = \infty$, cf.\ \ref{alg}\thetag{IIc}.

        If $r<\infty$ consider 
        \[
          A_{(r)}(t) := t^{-r} A(t).
        \]
        If $p<r$ or $A_{(r)}(0) = 0$, set $q(V) := 0$ and $V$ becomes a leaf of $T$. 
        Otherwise $A_{(r)}$ is $C^{p-r}$, and if $\mu_i$ is a choice of the eigenvalues for $A_{(r)}$, 
        then $\la_i(t)=t^r\mu_i(t)$ represent the eigenvalues of $t \mapsto A(t)$. Eigenvectors of $A_{(r)}$ are also 
        eigenvectors of $A$ (and vice versa). 
        Set $q(V):= r$, associate $V$ with $A_{(r)}$, and
        feed $A_{(r)}$ into \thetag{I}.
  \end{enumerate}
\end{algorithm}

We say that $T(A,0)$ is \emph{admissible} if the label $(d,q)$ of every leaf of $T$ satisfies
\begin{itemize}
  \item $q = 0 \Rightarrow d=1$,
\end{itemize}
i.e.,in Step~\thetag{2} we have $p \ge r$ and $A_{(r)}(0)\ne 0$.  
Evidently, admissibility is preserved by any of the reductions 
$A \leadsto A_i$, $A \leadsto A - \frac{1}{n} \on{Trace}(A) \I_n$, and $A \leadsto A_{(r)}$. 

Similarly we may consider $T(A,t)$ for any $t$.

\subsection{Associated integer \texorpdfstring{$\Th(A)$}{Theta(A)}} 
Let $A(t)$, $t \in I$, be a $C^p$-curve of normal matrices, $p \in \N \cup \{\infty\}$.
For each $t$ and each vertex $V$ in $T(A,t)$ we recursively define an integer $\Th_{t}(V)$ by setting
\begin{align*} 
    \Th_{t}(V) := 
      \begin{cases}
        q(V) & \text{ if $V$ is a leaf } \\
        q(V) + \max_W \Th_{t}(W) & \text{ otherwise }
      \end{cases}, 
\end{align*}
where $W$ ranges over all successors of $V$, and we set 
\begin{align} \label{eq:Th_alt}
    \Th_{t}(A) &:= \Th_{t}(\text{root of }T(A,t)) 
\end{align}
Furthermore, we set 
\begin{gather}
\Th(A) := \sup_{t \in I} \Th_{t}(A) \in \N \cup \{\infty\}. \label{Th_sup}
\end{gather}
It is easy to see that 
\[
  \Th(A) \le \ga(P_A) \le \Ga(P_A).
\]

\medskip

For later reference we also present an equivalent definition of $\Th_{t}(A)$ under the assumption 
that $T(A,t)$ is admissible:
\begin{itemize}
	\item If $A(t)$ has distinct eigenvalues and $A_i$, $1 \le i \le l$, denote the respective normal matrices introduced in 
		Algorithm~\ref{alg:n}\thetag{I}, 
		\begin{gather}
		\Th_{t}(A) := \max_i \Th_{t}(A_i). \label{Th_split}
		\end{gather}
	\item If all eigenvalues of $A(t)$ coincide, 
		replace $A$ by $A- \frac{1}{n} \on{Trace}(A) \I_n$ 
		(without changing $\Th_{t}(A)$). 
    Consider $r$ as in \eqref{eq:r}.
    \begin{itemize}
       \item If $r=\infty$ set $\Th_{t}(A) := 0$. 
       \item If $r<\infty$ set 
       \begin{align}
      \Th_{t}(A) &:= \Th_{t}(A_{(r)}) + r \label{Th_(r)}.
    \end{align}
    \end{itemize} 
\end{itemize}

\begin{remark*} 
We have:
  \begin{enumerate}
    \item If $A$ is admissible at $t$, then $A$ is $C^{\Th_{t}(A)}$ near $t$.
    \item If $A$ is $C^{1+\Th_{t}(A)}$ near $t$, then $A$ is admissible near $t$. 
  \end{enumerate}
\end{remark*}

\begin{theorem} \label{normal}
  Let $A(t)$, $t \in I$, be a $C^{p+\Th(A)}$-curve, $p \in \N \cup \{\infty\}$, of normal matrices such that $T(A,t)$ is admissible at any $t$ and 
  $E^{(\infty)}(A)=\emptyset$.
  Then the eigenvalues of $A$ can be parameterized by $C^{p+\Th(A)}$-functions, 
  the eigenvectors locally by $C^{p}$-functions and even globally if $p\ge 1$. 
\end{theorem}

\begin{demo}{Proof}
  We show the local assertions. Then the global statement for the eigenvalues follows from Lemma~\ref{glue}
  and Lemma~\ref{unique}. For the eigenvectors we may argue as in \cite[7.6]{AKLM98} or in \cite[5.11]{RainerN}; 
  here we need $p\ge 1$ 
  for the construction of the transformation function which involves linear ODEs. 
  Let $t$ be fixed.  
  
  \begin{claim*}[Local assertions]
    If $A$ is $C^{p+\Th(A)}$, then 
    the eigenvalues can be parameterized by $C^{p+\Th(A)}$-functions, the eigenvectors by $C^{p}$-functions, 
    locally near $t$.
  \end{claim*}
  
  By \eqref{Th_sup} it suffices to show the following claim.
  
  \begin{claim*}[1]
    If $A$ is $C^{p+\Th_{t}(A)}$, then 
    the eigenvalues can be parameterized by $C^{p+\Th_{t}(A)}$-functions, the eigenvectors by $C^{p}$-functions, 
    locally near $t$.
  \end{claim*}
  
  Without loss of generality assume that $0 \in I$ and $t=0$.
  We obtain Claim~\thetag{1} if we set $r=0$ and $A_{[0]}=A$ in the following claim.
  
  \begin{claim*}[2]
    Let $A_{[r]}$ be a $\tC{r}{p+\Th_{0}(A_{[r]})+r}$-curve, $r \in \N$, of normal matrices so that $T(A_{[r]},0)$ is admissible and 
    $0 \not\in E^{(\infty)}(A_{[r]})$.
    The eigenvalues of $A_{[r]}$ can be parameterized 
    by $\tC{r}{p+\Th_{0}(A_{[r]})+r}$-functions, the eigenvectors by $C^{p}$-functions.
  \end{claim*}

  Here $A_{[r]}$ should not be confused with $A_{(r)}$; 
  however, if $A \leadsto A_{(r)}$ as in Step~\thetag{IIa} of Algorithm~\ref{alg:n}, 
  then $A_{(r)}$ satisfies the assumptions of the claim.
  
  We prove Claim~\thetag{2} by induction on the size of the matrix.
  
  If all eigenvalues of $A_{[r]}(0)$ coincide, replace $A_{[r]}$ by $A_{[r]} - \frac{1}{n} \on{Trace}(A_{[r]}) \I_n$.
  By \ref{prop}\thetag{4}, this preserves the class $\tC{r}{p+\Th_{0}(A_{[r]})+r}$.
  Then all eigenvalues of $A_{[r]}(0)$ are equal to $0$, and $A_{[r]}(0)=0$.
  If $A_{[r]}=0$, then all eigenvalues vanish identically, the eigenvectors can be chosen constant, and we are done.
  If $A_{[r]} \ne 0$, we have $A_{[r]}(t) = t^s A_{(r+s)}(t)$ for $s \in \N_{>0}$, and $A_{(r+s)}$ is $\tC{r+s}{p+\Th_{0}(A_{[r]})+r}$
  with $A_{(r+s)}(0) \ne 0$, since $T(A_{[r]},0)$ is admissible. 
  Note that $T(A_{[r]},0) = T(A_{(r+s)},0)$, $\Th_0(A_{[r]})=\Th_0(A_{(r+s)}) + s$, by \eqref{Th_(r)},  
  $0 \not\in E^{(\infty)}(A_{(r+s)})$, 
  the eigenvalues of $A_{[r]}$ differ from those of $A_{(r+s)}$ by multiplication with $t^s$, and 
  eigenvectors of $A_{(r+s)}$ are eigenvectors of $A_{[r]}$ (and vice versa).  
  Replace $A_{[r]}$ by $A_{(r+s)}$.
  
  So after these initial steps we may assume that $A_{[r]}(0) \ne 0$ and $\on{Trace}(A_{[r]})=0$, thus 
  not all eigenvalues of $A_{[r]}(0)$ coincide.   
  By Step~\thetag{I} of Algorithm~\ref{alg:n} and by \ref{prop}\thetag{3}, the problem is reduced to  
  $\tC{r}{p+\Th_{0}(A_{[r]})+r}$-curves of normal matrices $A_{[r],j}$ of strictly smaller size, since 
  there exist real analytic mappings $\Ph_j$ defined near $A_{[r]}(0) \in \C^{n^2}$ 
  such that $A_{[r],j} = \Ph_j \o A_{[r]}$. Each $A_{[r],j}$ satisfies the assumptions of Claim~\thetag{2}, and, 
  by induction and \eqref{Th_split}, there exist 
  $\tC{r}{p+\Th_{0}(A_{[r]})+r}$-functions which parameterize the eigenvalues and $C^p$-functions which parameterize the eigenvectors 
  of each $A_{[r],j}$, respectively. This completes the proof of Claim~\thetag{2} and hence of the theorem.
\qed\end{demo}

\section{Definable normal matrices} \label{sec:defnormal}

We give a version of Theorem~\ref{normal} for the eigenvalues of definable curves of normal matrices;  
admissibility of the associated trees as well as
the assumption $E^{(\infty)}(A) = \emptyset$    
are not necessary. 
The eigenvectors, however, may not admit continuous parameterizations in this situation, see Example~\ref{ex:normal2}.  

Let $A(t)$, $t \in I$, be a definable $C^p$-curve of normal matrices.
We slightly modify Algorithm~\ref{alg:n}:
\begin{enumerate}[(IIa)]
  \item[(II)] By Lemma~\ref{deflem} we always have $p \ge r$.
  If $A_{(r)}(0)=0$, set $q(V) := r$ 
  and $V$ becomes a leaf of $T$. The eigenvalues of $A_{(r)}$ admit a $C^0$-parameterization by Lemma~\ref{controots}.   
\end{enumerate}  
The integer $\Th_{t}(A)$ is defined by \eqref{eq:Th_alt}.
Definability implies that the set of $t \in I$ such that $\Th_{t}(A)>0$ is finite and hence
$\Th(A)$  
is an integer.

\begin{theorem} \label{defnormal}
  Let $A(t)$, $t \in I$, be a definable $C^{p+\Th(A)}$-curve of normal matrices, $p \in \N \cup \{\infty\}$.
  Then the eigenvalues of $A$ can be parameterized by $C^{p+\Th(A)}$-functions.
\end{theorem}

Under these assumptions the eigenvectors of $A$ generally do not admit continuous parameterizations, 
see Example~\ref{ex:normal2}.

\begin{demo}{Proof}
  This follows by applying the arguments used in the proof of Theorem~\ref{maindef} to the proof of Theorem~\ref{normal}. 
\qed\end{demo}

\section{Examples} \label{sec:ex}

The examples in this sections show that the statements of Theorem~\ref{main}, \ref{maindef}, \ref{normal}, \ref{defnormal}, and 
Corollary~\ref{cor1}, except perhaps \ref{main}\thetag{3}, are optimal in the following sense:
\begin{itemize}
  \item The condition $E^{(\infty)}(P) = \emptyset$ is necessary unless $P$ is definable. Without that condition 
  a $C^\infty$-curve of hyperbolic polynomials need not be $C^{1,\al}$-solvable for any $\al>0$, see \cite{BBCP06}.
  \item Admissibility of the associated trees $T(P,t)$ is necessary unless $P$ is definable, see Example~\ref{ex:adm}.
  \item A (definable) $C^{p+\Ga_t(P)}$-curve of polynomials, where $T(P,t)$ is good (resp.\ bad) 
  is in general not $C^{p+\ga_t(P)+1}$ (resp.\ $C^{\ga_t(P)+1}$) solvable, see Example~\ref{ex:rad} and Example~\ref{ex:roots}.  
  \item The eigenvectors do not admit continuous parameterizations near $t$ if $T(A,t)$ is non-admissible or if 
  $t \in E^{(\infty)}(A)$, even if $A$ is definable, see Example~\ref{ex:normal2}.
  \item The eigenvectors of a (definable) $C^{p+\Th_t(A)}$-curve of normal matrices do in general not admit 
  $C^{p+1}$-parameterizations near $t$, see Example~\ref{ex:normal}. 
\end{itemize}
For $p \in \N$ consider the function $f_p$ defined by
\begin{equation*} \label{f_p}
f_p(t) := \left\{
\begin{array}{ll}
t^{p+1} & \text{if } t \ge 0\\ 0 & \text{if } t <0
\end{array}
\right.
\end{equation*}
which belongs to $C^{p}(\R) \setminus C^{p+1}(\R)$ and is definable. We shall make repeated use of $f_p$ in the 
following examples.

\begin{example} \label{ex:adm}
The function $f$ defined by 
\[
  f(t) = \begin{cases}
            t^{5} (\sin^2 \log t + 1) & t > 0 \\
            0 & t \le 0
         \end{cases}
\]
belongs to $C^4(\R) \setminus C^{5}(\R)$. Consider $P(t)(z) = z^2-f(t)$.   
We have $m_0(f) = 4$, that is $m=r=2$, but $T(P,0)$ is non-admissible, since $t^{-4} f(t)|_{t=0} = 0$.
In fact there is no $C^2$-function $g$ so that $g^2=f$, although $P_{(2)}$ is $C^0$-solvable.  
Indeed the roots  
\[
  \la_\pm(t) = \begin{cases}
            \pm t^2 \sqrt{t (\sin^2 \log t + 1)} & t > 0 \\
            0 & t \le 0
         \end{cases}
\]
are real analytic for $t\ne 0$. An easy computation shows that 
for $t_k = \exp(2\pi k)$, $k \in \Z$, we obtain 
\[
   \la_\pm''(t_k) \to \pm1 \quad \text{ as } \quad k \to -\infty,
\] 
and thus $\la_\pm$ cannot be $C^2$ at $t=0$.
\end{example}

\begin{example} \label{ex:rad}
Consider $f(t):= t^m (1+f_p(t))$. Then $m_0(f)=m$ and $f$ belongs to $C^{p+m} \setminus C^{p+m+1}$.
If $\frac{m}{n} \in \N$, then there is a $C^{p+\frac{m}{n}}$-function $g$ such that $g^n=f$, but $g$ cannot be in 
$C^{p+\frac{m}{n}+1}$.
If $\frac{m}{n} \not\in \N$, then $g$ can be chosen in $C^{\lfloor \frac{m}{n} \rfloor}$, but not in 
$C^{\lfloor \frac{m}{n} \rfloor+1}$. 
The following diagram illustrates $T(P,0)$ and $T(P^{\pm,h},0)$ for the case $\N \not\ni \frac{m}{n} = \frac{k}{h}$, 
where $h$ and $k$ are 
coprime.
  
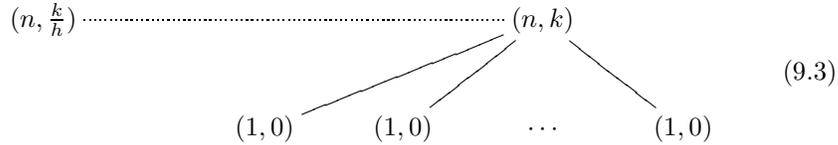
\begin{align} \label{diag:1}
  \begin{split}
    \xymatrix{
      (n,\frac{k}{h}) \ar@{..}[rrrr] &&&& (n,k) \ar@{-}[dll] \ar@{-}[dl] \ar@{-}[dr] & \\
      & & (1,0) & (1,0) & \cdots & (1,0)  
    }
  \end{split}
  \end{align} 
  \end{example}

\begin{example} \label{ex:roots}
Let $p,l \in \N$, $m=(n-1)l+s>0$ with $0 \le s < n-1$, and $p>n(l+1)$. Consider the $C^p$-curve of polynomials
\[
P_p(t)(z) = z^n - t^{m} z + f_p(t)\sum_{j=1}^n (-1)^j  z^{n-j}.
\]

Suppose $s=0$. 
Then $T(P_p,0)$ is good, $\Ga_0(P_p)=nl$, and $\ga_0(P_p)=l$. By Theorem~\ref{maindef}(1), $P_p$ admits 
$C^{p-(n-1)l}$-roots.
Suppose, for contradiction, that $P_p$ has $C^{p-(n-1)l+1}$-roots $\la_j$. 
Since $m_0(\la_j) \ge l$, we have $\la_j(t)=t^l \mu_j(t)$ for $C^{p-nl+1}$-functions $\mu_j$.
But then $f_p(t) = t^{nl} \mu_1(t) \cdots \mu_n(t)$ is $C^{p+1}$, by Lemma~\ref{deflem}, a contradiction. 
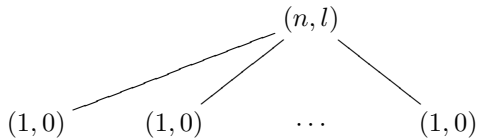
  \[
    \xymatrix{
      && (n,l) \ar@{-}[dll] \ar@{-}[dl] \ar@{-}[dr] & \\
       (1,0) & (1,0) & \cdots & (1,0)  
    } 
  \]

Suppose $s \ne 0$.
Then $T(P_p,0)$ is bad, $\Ga_0(P_p) = \lceil n (l+ \frac{s}{n-1}) \rceil  \le n(l+1)$, and $\ga_0(P_p)=l$. 
By Theorem~\ref{maindef}\thetag{2}, $P_p$ admits 
$C^{l}$-roots.
Suppose, for contradiction, that $P_p$ has $C^{l+1}$-roots $\la_j$.
Then $\la_j(t)=t^l \mu_j(t)$ for $C^{1}$-functions $\mu_j$ which represent the roots of $(P_p)_{(l)}$.
As $\mu_j(0)=0$ since $s>0$ (cf.\ Theorem~\ref{C^1}), 
the coefficient of $z$ in $P_p$ must have multiplicity $\ge (n-1)(l+1)$, a contradiction.  
The trees $T(P,0)$ and $T(P^{\pm,h},0)$ are given by diagram \eqref{diag:1} if we write $\frac{m}{n-1} = \frac{k}{h}$, 
where $h$ and $k$ are 
coprime.
\end{example}

\begin{example} \label{ex:normal2}
Consider the $C^p$-curve of symmetric matrices
  \[
  A(t) =t^{p+1} B(t) \quad \text{where}  \quad
  B(t)= 
  \begin{cases}
  \I_2 & t\ge 0 \\
  \begin{pmatrix}
     1& 1   \\
     1& -1
  \end{pmatrix} 
  & t<0 
  \end{cases}.
  \]
The tree $T(A,0)$ is not admissible, since $A_{(p)}(0)=0$, and $A$ does not admit a continuous parameterization of its eigenvectors.
The latter is true also for the definable $C^\infty$-curve of symmetric matrices $e^{-\frac{1}{t^2}} B(t)$. 
\end{example}

\begin{example} \label{ex:normal}
Consider the $C^{p+1}$-curve of symmetric matrices 
  \[
  A(t)= \begin{pmatrix}
    2 f_{p+1}(t) & t  \\
    t & 0
  \end{pmatrix},
  \quad t \in \R.
  \]
Its eigenvalues $f_{p+1}(t) \pm t \sqrt{1+f_{2p+1}(t)}$ are $C^{p+1}$ in accordance with Theorem~\ref{normal},  
as $\Th_0(A) = 1$. 
If $\binom{u}{v}$ denotes any $C^{p}$-eigenvector (which exists by Theorem~\ref{normal}), then 
\[
  u = (f_p \pm \sqrt{1+f_{2p+1}}) v.
\]
It follows that $u$ cannot be $C^{p+1}$ at $0$. Indeed, the Leibniz rule shows that the $(p+1)$st 
right handed derivative of $f_p v$ is non-zero because $v(0) \ne 0$, whereas $(f_p v)|_{t\le 0} = 0$. 
\end{example}

\def\cprime{$'$}
\providecommand{\bysame}{\leavevmode\hbox to3em{\hrulefill}\thinspace}
\providecommand{\MR}{\relax\ifhmode\unskip\space\fi MR }
\providecommand{\MRhref}[2]{%
  \href{http://www.ams.org/mathscinet-getitem?mr=#1}{#2}
}
\providecommand{\href}[2]{#2}

\end{document}